\documentclass[11pt]{amsart}
\usepackage{amsopn}
\usepackage{amssymb}
\usepackage{graphicx}
\usepackage[all]{xy}
\usepackage{amscd}
\usepackage{color}
\usepackage[colorlinks]{hyperref}
 \newtheorem{thm}{Theorem}[section]
 
 \newtheorem{lem}[thm]{Lemma}
 \newtheorem{prop}[thm]{Proposition}
 \newtheorem{rem}[thm]{Remark}
 \newtheorem{defn}[thm]{Definition}
 \newtheorem{exam}[thm]{Example}
 
 \theoremstyle{definition}
 \theoremstyle{remark}

\begin{document}

\title [Continuous K-g-frames]
{Continuous K-g-frames in Hilbert $C^*$-modules}

\author[J. Cheshmavar]{Jahangir Cheshmavar}

\author[J. Baradaran]{Javad Baradaran}

\author[A. Hosseinpour]{Asadollah Hosseinpour}

\address{Department of Mathematics, Payame Noor University (PNU),
P.O.Box: 19395-3697, Tehran, Iran}
\email{j$_{-}$cheshmavar@pnu.ac.ir}

\address{Department of Mathematics, Jahrom University, P.O.Box: 74135-111, Jahrom, Iran}
\email{baradaran@jahromu.ac.ir}

\address{Department of Mathematics, Payame Noor University (PNU), P.O.Box:
19395-3697, Tehran, Iran} \email{asadmath@student.pnu.ac.ir}

\thanks{}
\maketitle \dedicatory{}
\commby{}


\maketitle
\begin{abstract}
This study aims at combining the concepts of $g$-frame and
$K$-frame in a Hilbert $C^*$-module $U$ for an operator $K\in
End_{\mathcal{A}}^*(U)$, where ${\mathcal{A}}$ denotes a complex $C^{*}$-algebra and
$End_{\mathcal{A}}^*(U)$ denotes the set of all adjointable
$\mathcal{A}$-linear maps on $U$. As a result, the notion of a
$c$-$K$-$g$-frame for a Hilbert $C^{*}$-module is introduced and then characterizations on $c$-$K$-$g$-frames
 are provided. Finally, some results on a dual $c$-$K$-$g$-Bessel system in Hilbert $C^{*}$-modules are obtained.
\end{abstract}

 \maketitle
\section{\textbf{Introduction and preliminaries}}
Continuous frames were introduced in \cite{Ali, Kaiser}. A reason to study
frames in Hilbert $C^*$-modules is that there exist
differences between the structures of Hilbert spaces and Hilbert $C^*$-modules. For instance, in general, every bounded operator on a Hilbert
space has a unique adjoint while this  fact is not hold for  Hilbert
$C^*$-modules. Frank and Larson
\cite{M.Frank} extended the notion of a frame for an operator on a
Hilbert $C^*$-module. Henceforth, a number
generalizations of frames in Hilbert $C^*$-modules have attracted
more attention (for example, see
\cite{Alijani, Khosravi.2008, Xiao.Zeng} and the bibliographies
therein).
The concept of a $g$-frame for a Hilbert $C^*$-module
introduced in \cite{M.Frank, Khosravi.2008} and the
notion of a $K$-frame for a Hilbert $C^*$-module, for an operator $K
\in End_{\mathcal{A}}^*(\mathcal{H})$ defined in \cite{Najati}, where ${\mathcal{A}}$ is a complex $C^{*}$-algebra with the norm $||\;||_{\mathcal{A}}$ and $End_{\mathcal{A}}^*(\mathcal{H})$ denotes the set of all adjointable
$\mathcal{A}$-linear maps on a Hilbert $\mathcal{A}$-module $\mathcal{H}$. In what follows:
in section 2, the  purpose is to define the concept of a $c$-$K$-$g$-frame
for Hilbert $C^{*}$-modules and provide some results to the situation of $c$-$K$-$g$-frames in Hilbert $C^{*}$-modules (for example, Theorems \ref{thm1}, \ref{thm2}, \ref{thm3} ).
 In section 3, we  first introduce the dual of a $c$-$K$-$g$-frame in a Hilbert $C^{*}$-module and then we provide some characterizations of the dual $c$-$K$-$g$-Bessel system (Theorems \ref{thm4}, \ref{thm5}).
Throughout the work, we use $I$ as an indexing set and  $\mathbb{N}$, $\mathbb{C}$ as the set of natural numbers and complex numbers, respectively. We apply $U$ and $V$ for two Hilbert $C^*$-modules and
 $\{V_k\}_{k \in I}$ for a family of subspaces of $V$. For every $k \in I$,
  $End_{\mathcal{A}}^*(U, V_k)$ denotes
  the class of all adjointable $\mathcal{A}$-linear maps from
$U$ into $V_k$, also $End_{\mathcal{A}}^*(U)$ is abbreviated for $End_{\mathcal{A}}^*(U, \;U)$. The symbol $Ran(T)$ for the range of the operator $T$ is used. Finally, we assume
$(\mathcal{M}, \;\mu)$ is a measure space with positive measure
$\mu$. A complex linear
space $\mathcal{X}$ is called an (left) $\mathcal{A}$-module provided
that there exists a multiplication $\cdot: \mathcal{A}
\times \mathcal{X} \rightarrow \mathcal{X}$ having the following
properties:\\
 For all $f,g \in \mathcal{X}$, for all $a, b \in
\mathcal{A}$, and any $\lambda \in \mathbb{C}$,

\bigskip

\begin{itemize}
\item [{(i)}] $ a\cdot (f+g)=a\cdot f+a\cdot g$.~ \quad (ii)
$(a+b)\cdot f=a\cdot f+b\cdot f$~. \item [(iii)] $(a b)\cdot
f=a\cdot (b \cdot f)$~. \quad \quad \quad \,\ (iv) $\lambda(a
f)=(\lambda a)f=a(\lambda f)$ .~
\end{itemize}
\bigskip

A pre-Hilbert $\mathcal{A}$-module $\mathcal{X}$ is an (left)
$\mathcal{A}$-module equipped with an $\mathcal{A}$-valued inner
product $\langle \cdot, \cdot \rangle:\mathcal{X}\times
\mathcal{X}\longrightarrow \mathcal{A}$ satisfying:
\begin{itemize}
\item[(i)] $\langle f, \;f \rangle \geq 0, \;\forall f \in
\mathcal{X}$ and
$\langle f, \; f \rangle = 0$ if and only if $f=0$.~\\
\item[(ii)] $\langle f, \;g \rangle =\langle g, \;f \rangle^*, \;\forall
 f,\; g \in \mathcal{X}$~.\\
\item[(iii)] $\langle af+g, \;h \rangle= a \langle f, \;h \rangle+
\langle g, \;h \rangle, \; \forall a \in \mathcal{A}$, $\forall f,\; g, \;h
\in \mathcal{X}$.
\end{itemize}

So, it can be conclude that for any $a \in
\mathcal{A},$
\begin{eqnarray*}
\langle f, \;a g+h \rangle= \langle f, \;g \rangle a^*+
\langle f,\; h \rangle,\,\  \forall f,\; g,\; h \in \mathcal{X}.
\end{eqnarray*}
The norm on $\mathcal{X}$ is
defined by $\|f\|_{\mathcal{X}}=\|\langle f
f\rangle\|_{\mathcal{A}}^{1/2} \; for\; all f\in \mathcal{X}$. If $\mathcal{X}$ is a
Banach space with the norm $\|\cdot\|_{\mathcal{X}}$, then
$\mathcal{X}$ is said to be a Hilbert $\mathcal{A}$-module or Hilbert
$C^*$-module over $\mathcal{A}$. Since every complex Hilbert space is
a Hilbert $\mathbb{C}$-module,  Hilbert $C^*$-modules lay
between Hilbert spaces and Banach spaces. For a given $C^*$-algebra
$\mathcal{A}$, let $\bigoplus_{m\in \mathcal{M}}V_m$ be the
Hilbert $\mathcal{A}$-module is defined by
\begin{eqnarray*}
\bigoplus_{m\in \mathcal{M}}V_m:=\left\{\{g_m\}_{m\in \mathcal{M}}: g_m\in V_m,\,\
\|\int_{\mathcal{M}} |g_m|^2d\mu(m)\|<\infty \right\},
\end{eqnarray*}
 where for all $f = \{f_m\}_{m \in \mathcal{M}}$,  $g = \{g_m\}_{m \in
\mathcal{M}} \; \in  \bigoplus_{m\in \mathcal{M}}V_m$,
$\mathcal{A}$-valued inner product is given by
\begin{eqnarray*}
\langle f, \; g\rangle=\int_ {\mathcal{M}}\langle f_m, \;g_m\rangle
d\mu(m).
\end{eqnarray*}
We may refer the reader for more details  on $\bigoplus_{m\in \mathcal{M}}V_m$ to \cite{E.C. Lance}.\\

\begin{defn}
Let $U$ be a Hilbert $C^*$-module and $(\mathcal{M},\;
\mu)$ be a measure space with positive measure $\mu$. A continuous
frame, or simply a $c$-frame for $U$ with respect to $(\mathcal{M},
\mu)$ is a family $\{f_m\}_{m \in \mathcal{M}}$ of vectors in $U$ for which
\begin{itemize}
\item [(i)] For all $f \in U$, the mapping $m \mapsto \langle f, f_m \rangle$
is  measurable on $\mathcal{M}.$~
 \item [(ii)] There exist constants $A, B>0$ such that
\begin{eqnarray*}
 A\langle f, \; f \rangle \leq \int_{\mathcal{M}}
\langle f, \; f_m\rangle \langle f_m, \; f\rangle  d\mu(m)\leq B\langle
f, \; f\rangle, \; \forall f\in U.
\end{eqnarray*}
\end{itemize}
The constants $A, B$ are called $c$-frame bounds. The family $\{f_m\}_{m \in
\mathcal{M}}$ is called a $c$-Bessel system with Bessel bound $B$
if at least the right-hand of the above inequality is satisfied.
\end{defn}
\begin{defn}
A family $\{\Lambda_m\}_{m\in \mathcal{M}} \subset End_{\mathcal{A}}^*(U, V_m)$ is said to be a
continuous $g$-frame, or simply a $c$-$g$-frame for $U$ with respect to $\{V_m\}_{m \in
\mathcal{M}}$ if the following hold:
\begin{itemize}
\item [(i)] For all $f \in U$, the mapping $m \mapsto \langle f,\; \Lambda_mf \rangle$ is
measurable on $\mathcal{M}$.~

\item [(ii)] There exist constants $A, B>0$ such that
\begin{eqnarray*}
 A\langle f, \;f \rangle \leq \int_{\mathcal{M}}
\langle \Lambda_mf, \;\Lambda_mf\rangle d\mu(m)\leq B\langle f, \;
f\rangle, \; \forall f\in U.
\end{eqnarray*}
\end{itemize}
The constants $A, B$ are called $c$-$g$-frame bounds. If $A=B$, the family $\{\Lambda_m\}_{m\in \mathcal{M}}$ is said to be an $A$-tight $c$-$g$-frame for $U$.
A Parseval $c$-$g$-frame for $U$ with respect to $\{V_m\}_{m \in
\mathcal{M}}$  is an $A$-tight $c$-$g$-frame whenever $A$=1.
\end{defn}
Similar to the discrete case a continuous $g$-orthonormal
basis, or simply $c$-$g$-orthonormal basis for a family of adjointable operators on a Hilbert $C^{*}$-module
 is introduced (see \cite{Baradaran.2022}).
\begin{defn}
Let $\Theta=\{\Theta_{m}\}_{m\in \mathcal{M}}\subset End^{*}_{\mathcal{A}}(U, \; V_{m})$ be a family. \\
\begin{itemize}
\item [(1)] We call that $\Theta$ is $g$-complete if
 $\{f : \Theta_{m}f=0, \forall m\in
\mathcal{M}\}=\{0\}.$ \\
\item [(2)] We say  that $\Theta$ is a $c$-$g$-orthonormal basis for
$U$ with respect to
$\{V_m\}_{m \in \mathcal{M}}$ if it satisfies the following conditions:
\end{itemize}
\begin{itemize}
 \item [(i)] For all $f\in U$, the mapping $m\mapsto \langle f, \;\Theta_{m}f \rangle $ is measurable on
$\mathcal{M}$.~\\
  \item [(ii)] For almost all $n\in \mathcal{M}$,\\
   $\int_{\mathcal{M}}\langle \Theta^{*}_{m}g_{m}, \; \Theta^{*}_{n}g_{n}\rangle d\mu(m)=
   \langle g_{m},\; g_{n} \rangle, \; \forall g_{m}\in V_{m}, \;\forall g_{n}\in V_{n}.~$\\

 \item [(iii)] For any $f\in U$, $\int_{\mathcal{M}} \langle \Theta_{m}f, \; \Theta_{m}f \rangle d\mu(m)
 = \langle f,\; f \rangle$.
 \end{itemize}
\end{defn}
\begin{rem}
It is clear that any $c$-$g$-orthonormal basis is $g$-complete. In addition, it is a Parseval $c$-$g$-frame, i.e.,
 the $c$-$g$-frame operator $S$ is the identity operator.
\end{rem}
\begin{lem}\label{form.2}
\cite{Fang.Yao.2009} Let $U$, $V$ and $W$ be Hilbert
$\mathcal{A}-moduls$. Also, let $T^{'}\in End_{\mathcal{A}}^*(W,
V)$ and $T \in End_{\mathcal{A}}^*(U, V)$ with $\overline{Ran}(T^{*})$ orthogonally
complemented. The following statements are equivalent.
\begin{itemize}
\item [(i)] $T^{'}T^{'*}\leq \lambda TT^*$ for some $\lambda> 0.$
 \item [(ii)] There
exists $\mu>0$ such that $\|T^{'*}z\|\leq \mu\|T^*z\|, \; \forall z
\in V$.
 \item [(iii)] There exists a $D \in End_{\mathcal{A}}^*(W,
U)$ such that $T^{'}=TD$, i.e., $TX=T^{'}$ has a solution.
 \item
[(iv)] $R(T^{'})\subseteq R(T)$.
\end{itemize}
\end{lem}

\bigskip

\section{\textbf{c-K-g-frames in Hilbert $C^*$-Modules}}

In this section, the concept of a $c$-$K$-$g$-frame for
Hilbert $C^*$-modules is introduced and some results on $c$-$K$-$g$-frames are derived.

\begin{defn}\label{def1}
Suppose  that $( \mathcal{M}, \;\mu)$ is a measure space
with  positive measure $\mu$ and $K \in End_{\mathcal{A}}^*(U)$. We
say that $\{\Lambda_m\}_{m\in \mathcal{M}} \subset End_{\mathcal{A}}^*(U, \;V_m)$ is
a continuous $K$-$g$-frame, or simply a
$c$-$K$-$g$-frame for $U$ with respect to $\{V_m\}_{m \in
\mathcal{M}}$ if the following conditions hold:
\begin{itemize}
\item [(i)] For all $f \in U$, the mapping $m \mapsto \langle f,\; \Lambda_mf \rangle$ is
measurable on $\mathcal{M}$.~

\item [(ii)] There exist constants $A, B>0$ such that
\begin{eqnarray}\label{formula1}
A\langle K^*f, \;K^*f \rangle \leq
\int_{\mathcal{M}} \langle \Lambda_mf, \;\Lambda_mf\rangle
d\mu(m)\leq B\langle f, \;f\rangle, \; \forall f\in U.
\end{eqnarray}
\end{itemize}
The constants $A, B$ are called $c$-$K$-$g$-frame bounds. We call that
$\{\Lambda_m \}_{m \in \mathcal{M}}$ is a
tight  $c$-$K$-$g$-frame if there exists a constant $A>0$
such that
\begin{eqnarray*}
 \int_{\mathcal{M}} \langle \Lambda_mf, \;\Lambda_mf
\rangle d\mu(m)=A\langle K^*f, \;K^*f\rangle,\,\ \forall f \in U.
\end{eqnarray*}
A Parseval $c$-$K$-$g$-frame for $U$ with respect to
$\{V_m\}_{m\in \mathcal{M}}$ is a tight $c$-$K$-$g$-frame whenever $A=1$. We say the family $\{\Lambda_m\}_{m
\in \mathcal{M}}$ is a $c$-$g$-Bessel system for $U$ with respect to
$\{V_m\}_{m\in \mathcal{M}}$  with $c$-$g$-Bessel bound $B$ if at least the
right-hand inequality of (\ref{formula1}) is satisfied.
\end{defn}
\begin{rem}
Every $c$-$g$-frame for $U$ with respect to $\{V_m \}_{m \in \mathcal{M}}$ is a $c$-$K$-$g$-frame for $U$. In fact, if
$\{\Lambda_m\}_{m \in \mathcal{M}}$ is a $c$-$g$-frame for $U$  with respect to $\{V_m \}_{m \in \mathcal{M}}$,
 then by the definition there exist constants $A, B> 0$ such that
\begin{eqnarray*}
A\langle f, \;f \rangle \leq \int_{\mathcal{M}} \langle \Lambda_mf, \;
\Lambda_mf\rangle d\mu(m)\leq B\langle f, \; f\rangle,\; \; \forall f \in U.
\end{eqnarray*}
On the other hand, using $\langle K^*f, \; K^*f\rangle \leq \|K\|^2
\langle f, \; f \rangle,\;  \forall f \in U$, we obtain
\begin{eqnarray*}
A\|K\|^{-2}\langle K^*f, \; K^*f \rangle \leq \int_{\mathcal{M}}
\langle \Lambda_mf, \; \Lambda_mf\rangle d\mu(m)\leq B\langle f, \;
f\rangle,\,\ \forall f \in U.
\end{eqnarray*}
Therefore, $\{\Lambda_m \}_{m \in \mathcal{M}}$ is a $c$-$K$-$g$-frame for
$U$ with respect to $\{V_m\}_{m \in \mathcal{M}}$.
\end{rem}
\begin{exam}
Let $(\mathcal{M}, \; \mu)$ be a measure space with
$\mathcal{M}=[0,\; 1]$ and the Lebesgue measure $\mu$. Assume that
\begin{eqnarray*}
U= \left \{
\begin{bmatrix}
a & 0 & \\
0 & b & \\
\end{bmatrix}, \; \;  a, b \in \mathbb{C}   \right \}.
\end{eqnarray*}
By a simple computation, $U$ is a unital $C^*$-algebra and also,
$U$ is a Hilbert $C^{*}$-module over itself with the following
inner product,
$$ \langle ., \; . \rangle : U\times U \rightarrow U,\; \; \langle A, \; B \rangle=A(\overline{B})^{t},$$
where $(\bar{B})^{t}$ denotes the conjugate transpose of  the matrix $B$. For
any $m \in \mathcal{M}$, we take $V_m=U$. Consider the family
$\{f_m\}_{m \in \mathcal{M}}$ as
\begin{eqnarray*}
f_{m}=
\begin{bmatrix}
\frac{m}{2} & 0 & \\
0 & m+1 & \\
\end{bmatrix},
\end{eqnarray*}
and the family $\{\Lambda_{m}\}_{m \in \mathcal{M}}$ as
$$\Lambda_{m}f:= \langle f, \; f_m \rangle,\,\ \forall f \in U.$$
Then, for any $f\in U$, we obtain
\begin{eqnarray*}
\int_{\mathcal{M}} \langle \Lambda_mf, \; \Lambda_mf \rangle
d\mu(m)&=&
\int_{\mathcal{M}} \left\langle f, \; f_{m} \right\rangle  \left\langle f_{m},\; f \right\rangle d\mu(m) \\
&=& \int_{[0,1]} \left\langle
\begin{bmatrix}
a &0& \\
0 &b& \\
\end{bmatrix}
,
\begin{bmatrix}
\frac{m}{2} & 0 & \\
0 & m+1 & \\
\end{bmatrix}
\right\rangle\\
&\times& \left\langle
\begin{bmatrix}
\frac{m}{2} & 0 & \\
0 & m+1 & \\
\end{bmatrix}
,
\begin{bmatrix}
a &0& \\
0 &b& \\
\end{bmatrix}
\right\rangle d\mu(m)\\ &=& \int_{[0,1]} \begin{bmatrix}
\frac{m}{2}a &0& \\
0 &(m+1)b& \\
\end{bmatrix}
\begin{bmatrix}
\frac{m}{2}\overline{a} &0& \\
0 &(m+1)\overline{b}& \\
\end{bmatrix}
d\mu(m)\\
&=& \int_{[0,1]} \begin{bmatrix}
\frac{m^2}{4} &0& \\
0 &(m+1)^2& \\
\end{bmatrix}
\begin{bmatrix}
|a|^2 &0& \\
0 &|b|^2& \\
\end{bmatrix}
d\mu(m)\\
&=& \begin{bmatrix}
|a|^2 &0& \\
0 &|b|^2& \\
\end{bmatrix}
\int_{[0,1]} \begin{bmatrix}
\frac{m^2}{4} &0& \\
0 &(m+1)^2& \\
\end{bmatrix}
d\mu(m)\\
&=& \begin{bmatrix}
\frac{1}{12} &0& \\
0 &\frac{8}{3}& \\
\end{bmatrix}
\begin{bmatrix}
|a|^2 &0& \\
0 &|b|^2& \\
\end{bmatrix}
\end{eqnarray*}
Hence, for every $f \in U$, we can conclude that
\begin{eqnarray*}
\frac{1}{12} \langle f,f \rangle \leq \int_{\mathcal{M}} \langle
\Lambda_mf, \; \Lambda_mf \rangle d\mu(m)\leq \frac{8}{3} \langle
f, f \rangle.
\end{eqnarray*}

On the other hand, for any $K \in End_{\mathcal{A}}^*(U)$, we get
$$\langle K^*f, K^*f \rangle \leq \|K\|^2 \langle
f, f \rangle, \; \; \forall f \in U.$$  Now, it follows from two the above inequalities
\begin{eqnarray*}
\frac{1}{12} \|K\|^{-2} \langle K^*f, K^*f \rangle \leq
\int_{\mathcal{M}} \langle \Lambda_mf, \; \Lambda_mf \rangle
d\mu(m)\leq \frac{8}{3} \langle f, f \rangle, \,\ \forall f \in U.
\end{eqnarray*}
That is, $\{\Lambda_{m}\}_{m \in \mathcal{M}}$ is a
$c$-$K$-$g$-frame for $U$ with respect to $\{V_m\}_{m\in
\mathcal{M}}$ and with the bounds $\frac{1}{12}\|K\|^{-2}$,
$\frac{8}{3}$.
\end{exam}

\begin{defn} \cite{Baradaran.2022}
Suppose  $K \in End_{\mathcal{A}}^*(U)$ and $\{\Lambda_m\}_{m\in \mathcal{M}}\subset End_{\mathcal{A}}^*(U,\; V_{m})$ is a $c$-$K$-$g$-frame for $U$ with respect to $\{V_m\}_{m \in \mathcal{M}}$. We define the linear operator $S$ as follows:
$$ S: U\longrightarrow U, \; \; \;  Sf=\int_{\mathcal{M}}\Lambda_{m}^{*}\Lambda_{m}f d\mu(m), \; \; \forall\; f \in U.$$
The operator $S$ is called the $c$-$g$-frame operator of $\{\Lambda_m\}_{m\in \mathcal{M}}$.
\end{defn}
According to Proposition $3.2$ in \cite{Rashidi.2011}, $S$ is a positive, bounded and self-adjoint operator. Furthermore, we can conclude that the operator $S$ satisfies:
\begin{eqnarray*}
\langle Sf, \;g\rangle = \int_{\mathcal{M}} \langle \Lambda_m^{*}
\Lambda_mf , \; g \rangle   d\mu(m), \,\ \forall f,\; g \in U.
\end{eqnarray*}
\begin{rem}
Recall that, in general, similar to $K$-frames, the $c$-$g$-frame operator $S$ of a $c$-$K$-$g$-frame for a Hilbert $C^{*}$-module $U$
is not invertible on $U$. However, if $K$ has closed range, then $S$ on $R(K)$ is a positive, self-adjoint, and invertible operator
(see \cite{Gavruta.2012, Xiao.Zhu}).
\end{rem}

\begin{lem} \cite{Baradaran.2022} \label{form.22}
Let $\Theta=\{\Theta_{m}\}_{m\in \mathcal{M}} \subset End^{*}_{\mathcal{A}}(U, \; V_{m})$ be
a $c$-$g$-orthonormal basis for $U$ with respect to $\{V_{m}\}_{m\in \mathcal{M}}$. Suppose
$\Lambda=\{\Lambda_{m}\}_{m\in \mathcal{M}}\subset End^{*}_{\mathcal{A}}(U, \; V_{m})$ is a family such that
the mapping $m \mapsto \langle f,\; \Lambda_{m}f \rangle$ is measurable
for any $f \in U.$ Then,
$\Lambda$ is a $c$-$g$-Bessel
system for $U$ with respect to $\{V_{m}\}_{m\in \mathcal{M}}$ if and only if there exits a unique bounded operator $R\in End^{*}_{\mathcal{A}}(U)$
with  $\Lambda_{m}=\Theta_{m}R^{*}, \; \forall m\in \mathcal{M}.$
\end{lem}
Assume that $\{\Lambda_m\}_{m \in \mathcal{M}} \subset End^{*}_{\mathcal{A}}(U, V_{m})$ is a $c$-$g$-frame for
$U$ with respect to $\{V_m\}_{m \in \mathcal{M}}$  with the $c$-$g$-frame operator $S$.
Similar to the usual frames, the family $\{\Lambda_{m}S^{-\frac{1}{2}}\}_{m \in \mathcal {M}}$
is a Parseval $c$-$g$-frame for $U$ with respect to $\{V_m\}_{m \in \mathcal{M}}$.
For a given $K\in End^{*}_{\mathcal{A}}(U)$, we denote $T_{K}(U)$ as the class of all tight $c$-$K$-$g$-frames
for a Hilbert $C^{*}$-module $U$.
\begin{thm}
Suppose that $K_{1}, K_{2}\in End^{*}_{\mathcal{A}}(U)$ and $U$ is a finitely or
countably generated Hilbert $C^{*}$-module.
Then, $T_{K_1}(U) \subseteq T_{K_2}(U)$ if and only if
$K_{1}K_{1}^{*}=\lambda K_{2}K_{2}^{*}$ for some $\lambda>0$.
\end{thm}
\begin{proof}
Let $\{\Lambda_m\}_{m \in \mathcal{M}}\subset End_{\mathcal{A}}^*(U, V_m)$ be a
 $c$-$g$-frame for $U$ with respect to
$\{V_{m}\}_{m \in \mathcal{M}}$ with $c$-$g$-frame operator $S$. Since
 $\{\Lambda_{m}S^{\frac{-1}{2}}\}_{m \in \mathcal{M}}$ is a Parseval $c$-$g$-frame for $U$,
we have
\begin{eqnarray*}
 \int_{\mathcal{M}} \langle \Lambda_mS^{-\frac{1}{2}}f, \;
\Lambda_mS^{-\frac{1}{2}}f \rangle d\mu(m)=\langle f,\; f \rangle, \; \forall f\in U.
\end{eqnarray*}
Hence, for any $f\in U$, we get
\begin{eqnarray}\label{formula2}
 \int_{\mathcal{M}} \langle \Lambda_mS^{-\frac{1}{2}}K^{*}_{1}f, \;
\Lambda_mS^{-\frac{1}{2}}K^{*}_{1}f \rangle d\mu(m)=\langle K^{*}_{1}f,\; K^{*}_{1}f \rangle.
\end{eqnarray}
Therefore, the family $\{\Lambda_{m}S^{-\frac{1}{2}}K^{*}_{1}\}_{m \in \mathcal{M}}$ is
a tight $c$-$K_{1}$-$g$-frame for $U$.\\
Firstly, let $T_{K_1}(U)\subseteq T_{K_2}(U)$,
 so $\{\Lambda_{m}S^{-\frac{1}{2}}K_{1}^{*}\}_{m \in \mathcal{M}}$ is a tight $c$-$K_{2}$-$g$-frame for $U$.
Thus, for any $f\in U$ and for some $\lambda>0$, we obtain
\begin{eqnarray}\label{formula3}
 \int_{\mathcal{M}} \langle \Lambda_m S^{-\frac{1}{2}}K_{1}^{*}f, \;
\Lambda_m S^{-\frac{1}{2}}K_{1}^{*}f \rangle d\mu(m)=\lambda \langle K^{*}_{2}f,\; K^{*}_{2}f \rangle.
\end{eqnarray}

 It follows from the relations (\ref{formula2}), (\ref{formula3})
 $$\langle K_{1}K^{*}_{1}f,\; f \rangle=\langle \lambda K_{2}K^{*}_{2}f,\; f \rangle, \; \forall f \in U.$$
Also, using the Polarization formula, we get
$$\langle K_{1}K^{*}_{1}f,\; g \rangle=\langle \lambda K_{2}K^{*}_{2}f,\; f \rangle, \; \forall f ,g \in U,$$
 i.e., $K_{1}K_{1}^{*}=\lambda K_{2}K_{2}^{*}.$\\
 In contrast, let $K_{1}K_{1}^{*}=\lambda K_{2}K_{2}^{*}$ for some $\lambda>0$. Let
 $\{\Lambda_m\}_{m \in \mathcal{M}}\subset End_{\mathcal{A}}^*(U, V_m)$ be a tight
  $c$-$K_{1}$-$g$-frame for $U$ with respect to $\{V_{m}\}_{m \in \mathcal{M}}$.
 Then, for some $A_{1}>0$, we have
 $$\int_{\mathcal{M}} \langle \Lambda_mf,\;
\Lambda_mf \rangle d\mu(m)= A_{1}\langle K^{*}_{1}f,\; K^{*}_{1}f \rangle , \; \forall f\in U.$$
 Thus, for each $f\in U$, we obtain $$\int_{\mathcal{M}} \langle \Lambda_mf, \; \Lambda_mf \rangle d\mu(m)
 =A_{1}\langle K^{*}_{1}f,\; K^{*}_{1}f \rangle
 =A_{1}\lambda \langle K^{*}_{2}f,\; K^{*}_{2}f \rangle.$$
Hence, $\{\Lambda_m\}_{m\in \mathcal{M}}$ is a tight $c$-$K_{2}$-$g$-frame for $U$,
so $T_{K_{1}}(U)\subseteq T_{K_{2}}(U).$
\end{proof}
\begin{prop}\label{prop2}
Every adjointable operator $K$ on a unital Hilbert $C^*$-module $U$
admits a $c$-$K$-$g$-frame for $U$.
\end{prop}
\begin{proof}
Let $\{f_m\}_{m \in \mathcal{M}}$ be a $c$-frame for $U$ with bounds $A, B$. For every $m \in \mathcal{M}$,
we take $V_m=\mathcal{A}$, and define adjointable operator $\Lambda_{m}$ as follows:
$$\Lambda_m: U \rightarrow V_m ,\;  \; \Lambda_m f=\langle f,\; f_m\rangle, \; \forall f\in U.$$
Using the following relation:
$$ \langle f,\;  f_{m}\rangle \langle f_{m},\; f\rangle 1_{\mathcal{A}}=\langle f,\; f_m\rangle \langle f_{m},
\; f\rangle \langle 1,\; 1 \rangle= \langle \langle f, \; f_{m}\rangle, \;  \langle f, \; f_{m}\rangle \rangle,$$
and the definition of  $c$-frame, we obtain
\begin{eqnarray*}
A\langle f, \; f\rangle \leq \int_{\mathcal{M}} \langle \Lambda_m f, \;
\Lambda_m f \rangle d\mu(m) \leq B \langle f, \; f\rangle,\; \forall f\in U.
\end{eqnarray*}
If for every $m \in \mathcal{M}$, we put $\Gamma_m=\Lambda_mK^*$, then for any $f\in U$, we get
\begin{eqnarray*}
A\langle K^*f, \; K^*f\rangle &\leq& \int_{\mathcal{M}} \langle
\Lambda_m K^*f, \; \Lambda_m K^*f \rangle d\mu (m)\\
&=& \int_{\mathcal{M}} \langle \Gamma_mf, \; \Gamma_mf \rangle
d\mu(m)\\
&\leq& B \langle K^{*} f, \; K^{*}f \rangle \\
&\leq& B \|K\|^2\langle f, \; f\rangle.
\end{eqnarray*}
By the definition of $ \{\Lambda_m\}_{m \in \mathcal{M}}$, since for every $f \in U$, the mapping $m \mapsto \langle f,\; \Gamma_mf \rangle$ is
measurable on $\mathcal{M}$,~  $\{\Gamma_m\}_{m \in \mathcal{M}}$ is a $c$-$K$-$g$-frame
for $U$ with respect to $\{V_m\}_{m \in \mathcal{M}}$.
\end{proof}
\begin{thm}
Assume $K_1,  K_2 \in End_{\mathcal{A}}^* (U)$ with
$R(K^{*}_1)\perp R(K^{*}_2)$. If $\{\Lambda_{m}\}_{m \in
\mathcal{M}}\subset End_{\mathcal{A}}^*(U, V_m)$ is a
$c$-$K_1$-$g$-frame as well as a $c$-$K_2$-$g$-frame for $U$ with
respect to $\{V_m\}_{m \in \mathcal{M}}$ and $\alpha , \beta$ are
scalars, then $\{\Lambda_{m}\}_{m \in \mathcal{M}}$ is a
$c$-$(\alpha K_1+\beta K_2 )$-$g$-frame and a $c$-$K_1K_2$-$g$-frame
for $U$ with respect to $\{V_{m}\}_{m \in \mathcal{M}}$.
\end{thm}
\begin{proof}
Since $\{\Lambda_{m}\}_{m \in \mathcal{M}}$ is a $c$-$K_1$-$g$-frame
and a $c$-$K_2$-$g$-frame for $U$ with respect to
$\{V_{m}\}_{m \in \mathcal{M}}$, there exist positive
constants $A_n$ and $B_n (n = 1,2)$ such that for any $f \in U$,
we have
\begin{eqnarray*}
A_n \langle K_n^* f, \; K_n^* f \rangle \leq \int_{\mathcal{M}}
\langle \Lambda_mf,\; \Lambda_mf\rangle d\mu(m)
\leq B_n\langle f, \; f\rangle.
\end{eqnarray*}
Now, for each $f \in U$, we can write
\begin{eqnarray*}
\langle (\alpha K_1 + \beta K_2)^* f,\; (\alpha K_1 + \beta K_2)^* f
\rangle = \langle \overline{\alpha} K_1^*f + \overline{\beta}
K_2^*f,\; \overline{\alpha} K_1^*f + \overline{\beta} K_2^* f
\rangle\\
=|\alpha|^2 \langle K_1^*f,\; K_1^*f \rangle + \overline{\alpha}
\beta \langle K_1^*f, \; K_2^*f \rangle+ \alpha \overline{\beta}
\langle K_2^*f,\; K_1^*f \rangle+|\beta|^2 \langle K_2^*f,\;K_2^*f
\rangle.
\end{eqnarray*}
Because $R(K^{*}_1)\perp R(K^{*}_2)$, two terms in the middle of the
above equality are zero, the above equality continues as way:
\begin{eqnarray*}
 =|\alpha|^2 \langle K_1^*f, \; K_1^*f \rangle + |\beta|^2 \langle
K_2^*f, \; K_2^*f \rangle.
\end{eqnarray*}
Therefore, for any $f \in U$, we obtain
\begin{eqnarray*}
\frac{A_1A_2}{2(|\alpha|^2A_2+ |\beta|^2A_1)} \langle (\alpha K_1
+ \beta K_2)^* f, \; (\alpha K_1 + \beta K_2)^* f \rangle\\
 =\frac{A_1A_2|\alpha|^2}{2(|\alpha|^2A_2+ |\beta|^2A_1)} \langle
K_1^*f, \; K_1^*f \rangle + \frac{A_1A_2|\beta|^2}{2(|\alpha|^2A_2+
|\beta|^2A_1)} \langle K_2^*f, \; K_2^*f \rangle\\
\leq \frac{1}{2}\int_{\mathcal{M}} \langle \Lambda_mf, \;
\Lambda_mf\rangle d\mu(m)+
\frac{1}{2}\int_{\mathcal{M}} \langle \Lambda_mf,\; \Lambda_mf\rangle d\mu(m)\\
\leq (\frac{B_1+B_2}{2}) \langle f,\; f\rangle.
\end{eqnarray*}
Thus, $\{\Lambda_{m}\}_{m \in \mathcal{M}}$ is a $c-(\alpha
K_1 + \beta K_2)$-$g$-frame for $U$ with respect to $\{V_{m}\}_{m
\in \mathcal{M}}$. To do the second part, on the one hand, for each $f \in U$, we have
\begin{eqnarray*}
\langle (K_1K_2)^*f,\; (K_1K_2)^*f \rangle &=& \langle
K_2^*K_1^*f, \;K_2^*K_1^*f \rangle\\
 &\leq& \Arrowvert K_2^*
\Arrowvert^2 \langle K_1^*f, \;K_1^*f \rangle.
\end{eqnarray*}
On the other hand,  since $\{\Lambda_{m}\}_{m \in \mathcal{M}}$ is a
$c$-$K_1$-$g$-frame for $U$ with respect to $\{V_{m}\}_{m \in
\mathcal{M}}$, for any $f \in U$, we get
\begin{eqnarray*}
\frac{A_1}{\Arrowvert K_2^* \Arrowvert^2} \langle
(K_1K_2)^*f,\; (K_1K_2)^*f \rangle \leq \int_{\mathcal{M}} \langle
\Lambda_mf, \; \Lambda_mf\rangle d\mu \leq B_1 \langle f,\; f\rangle.
\end{eqnarray*}
 Hence, $\{\Lambda_{m}\}_{m \in \mathcal{M}}$ is
a $c$-$K_1K_2$-$g$-frame for $U$ with respect to $\{V_{m}\}_{m \in
\mathcal{M}}$.
\end{proof}

\begin{thm}\label{thm1}
Let $\{\Lambda_m\}_{m \in \mathcal{M}}$ be a $c$-$g$-Bessel system
for $U$ with respect to $\{V_m\}_{m \in \mathcal{M}}$ and with $c$-$g$-Bessel bound
$B$. Then, for all $f=\{f_{m}\}_{m\in \mathcal{M}}\in \bigoplus_{m\in \mathcal{M}}V_m$,
the mapping $T: \bigoplus_{m\in \mathcal{M}}V_m \rightarrow U$
defined by
\begin{eqnarray*}
\langle Tf, \; g \rangle=\int_{\mathcal{M}} \langle \Lambda_m^*f_m, \;
 g \rangle d\mu(m),\; \; \forall g \in U,
\end{eqnarray*}
is a bounded linear operator with $\|T\| \leq \sqrt{B}$. Moreover, the adjoint operator is  $T^{*}g=\{\Lambda_{m}g\}_{i \in I}$ for all $g \in U$ with the property $S=TT^{*}.$
\end{thm}
The operator $T$ is called the synthesis operator of $c$-$g$-Bessel system $\{\Lambda_m\}_{m \in \mathcal{M}}$.
\begin{proof}
By a simple calculation, one can see that the mapping $T$ is linear and adjointable.
For every $f=\{f_{m}\}_{m\in \mathcal{M}}\in \bigoplus_{m\in \mathcal{M}}V_m$, we obtain
\begin{eqnarray*}
||Tf||&=&sup_{||g||\leq 1}||\langle Tf,\; g \rangle||= sup_{||g||\leq 1} \left|\left|\int_{\mathcal{M}} \langle \Lambda_m^*f_m, \;
 g \rangle d\mu(m) \right|\right| \\
 &=&sup_{||g||\leq 1} \left|\left|\int_{\mathcal{M}} \langle f_m, \; \Lambda_mg \rangle d\mu(m) \right|\right| \\
 &\leq& sup_{||g||\leq 1} \left (\int_{\mathcal{M}} ||\Lambda_mg||^{2}d\mu(m) \right)^{\frac{1}{2}} \left(\int_{\mathcal{M}} ||f_m||^{2}d\mu(m)\right)^{\frac{1}{2}}\\
 &\leq&\sqrt{B}||f||.
\end{eqnarray*}
So, $||T||\leq \sqrt{B}.$
Also, for any $f=\{f_{m}\}_{m\in \mathcal{M}}\in \bigoplus_{m\in \mathcal{M}}V_m$  and for any $g \in U$,
we have
$$\langle f, \; \{\Lambda_{m}g\}_{m\in \mathcal{M}} \rangle
=\int_{\mathcal{M}}\langle \Lambda_{m}^{*}f_{m}, \; g \rangle  d\mu(m)=\langle Tf, \; g \rangle.$$
It follows from this the adjoint operator is
$T^{*}(g)=\{\Lambda_{m}g\}_{m\in \mathcal{M}}, \;\; \forall \; g\in U.$
Furthermore, for all $h, \; g \in U$, we get
$$\langle TT^{*}g, \; h \rangle=\int_{\mathcal{M}}\langle \Lambda_{m}^{*}\Lambda_{m}g, \; h \rangle d\mu(m)
=\langle Sg, \; h \rangle,\; i.e., \; S=TT^{*},$$ where $S$ is the $c$-$g$-frame operator of
$\{\Lambda_{m}\}_{m\in \mathcal{M}}.$
\end{proof}
\begin{rem}\label{rem10}
If we define the operator $R: \bigoplus_{m\in \mathcal{M}}V_m \rightarrow U$ by
$$Rf=\int_{\mathcal{M}} \Lambda^{*}_{m}f_{m} d\mu(m),
 \; \; \forall f=\{f_{m}\}_{m\in \mathcal{M}}\in \bigoplus_{m\in \mathcal{M}}V_m,$$
then for each $g\in U$, we conclude that  $R^{*}g=\{\Lambda_{m}g\}_{m\in  \mathcal{M}}$.  Moreover, for every  $h\in U$, we have
 $\langle RR^{*}h, \; g \rangle= \langle Sh, \; g \rangle$,i.e., $S=RR^{*}$.
 Therefore, it follows that $T=R$.
\end{rem}
\begin{rem}
Applying Theorem \ref{thm1}, the unique operator $R$ in
Lemma \ref{form.22} is exact the synthesis operator of $c$-$g$-Bessel
system $\Lambda=\{\Lambda_{m}\}_{m\in \mathcal{M}}$. In fact, the operator $R$  and the operator
$T$  in Theorem \ref{thm1} have the same values,
though those are defined in different ways.
\end{rem}

The next two results  provide  characterizations of $c$-$K$-$g$-frames.
\begin{thm}\label{thm2}
Let $K \in End_{\mathcal{A}}^*(U)$ and $\{\Lambda_m\}_{m\in \mathcal{M}}\subset End_{\mathcal{A}}^*(U, V_m)$
such that for any $f \in U$, the mapping $m \mapsto \langle f, \Lambda_mf \rangle$
is  measurable on ${\mathcal{M}}.$  Suppose the operator $T^{*}: U \mapsto \bigoplus_{m\in \mathcal{M}}V_m$  defined by $T^{*}f=\{\Lambda_{m}f\}_{m \in \mathcal{M}}$ and
$\overline{Ran}(T^{*})$ is orthogonally complemented. Then, $\{\Lambda_m\}_{m\in \mathcal{M}}$ is a $c$-$K$-$g$-frame for $U$ with respect to $\{V_m\}_{m \in \mathcal{M}}$
if and only if there are constants $A, B > 0$ such that
\begin{eqnarray}\label{formula4}
A||K^{*}f||^{2}\leq \left|\left|\int_{\mathcal{M}}{\langle \Lambda_{m}f,\; \Lambda_{m}f \rangle}d\mu(m) \right|\right|\leq B||f||^{2},  \; \; \forall f \in U.
\end{eqnarray}
\end{thm}
\begin{proof}
Let $\{ \Lambda_{m}\}_{m \in \mathcal{M}}$ be a
$c$-$K$-$g$-frame for $U$ with respect to $\{V_m\}_{m \in
\mathcal{M}}$. It is clear that $\{ \Lambda_{m}\}_{m \in \mathcal{M}}$  satisfies the relation (\ref{formula4}).\\
Conversely, assume that (\ref{formula4}) holds. The left-hand inequality of (\ref{formula4}) implies that $||K^{*}f||^{2}\leq \frac{1}{A}||T^{*}f||^{2}, \; \forall f \in U.$
Since $T \in End_{\mathcal{A}}^*(\bigoplus_{m\in \mathcal{M}}V_m, U)$, using Lemma \ref{form.2}, there exists a constant $\lambda> 0$ such that $KK^{*}\leq\lambda TT^{*}.$
It follows from  that
\begin{eqnarray}\label{formula5}
\frac{1}{\lambda}\langle K^{*}f,\; K^{*}f \rangle \leq \langle T^{*}f,\;  T^{*}f \rangle=\int_{\mathcal{M}}{\langle \Lambda_{m}f,\; \Lambda_{m}f \rangle}d\mu(m), \; \forall f \in U.
\end{eqnarray}
Now, we show that $\{ \Lambda_{m}\}_{m \in \mathcal{M}}$ is a
$c$-$g$-Bessel system for $U.$ To this, for any $\{g_{m}\}_{m \in \mathcal{M}} \in \bigoplus_{m\in \mathcal{M}}V_m $ and any $I\subset \mathcal{M}$, the right-hand inequality of
(\ref{formula4}) results that
\begin{eqnarray*}
\left|\left|\int_{I}\Lambda^{*}_{m}g_{m}d\mu(m) \right|\right|^{2}&=&sup_{f \in U, ||f||=1}\left|\left|\langle \int_{I}\Lambda^{*}_{m}g_{m}d\mu(m),\; f \rangle \right|\right|^{2}\\
&=&sup_{f \in U, ||f||=1}\left|\left|\int_{I} \langle \Lambda^{*}_{m}g_{m},\; f \rangle d\mu(m) \right|\right|^{2} \\
&=&sup_{f \in U, ||f||=1}\left|\left|\int_{I} \langle g_{m},\; \Lambda_{m}f \rangle d\mu(m) \right|\right|^{2}\\
&\leq& sup_{f\in U, ||f||=1}\left|\left|\int_{I} \langle g_{m},\;
g_{m} \rangle d\mu(m) \right|\right|\\
&\times& \left|\left|\int_{I} \langle \Lambda_{m}f,\; \Lambda_{m}f \rangle d\mu(m) \right|\right|\\
&\leq& B \left|\left|\int_{I} |g_{m}|^{2} d\mu(m)
\right|\right|<\infty.
\end{eqnarray*}
So, $\int_{\mathcal{M}}\Lambda^{*}_{m}g_{m}d\mu(m)$ converges in $H$. Hence, we have
\begin{eqnarray*}
\langle T^{*}f,\; \{g_{m}\}_{m \in \mathcal{M}} \rangle&=& \int_{\mathcal{M}} \langle \Lambda_{m}f,\; g_{m} \rangle d\mu(m)\\
&=&\langle f,\; \int_{\mathcal{M}}\Lambda_{m}^{*}g_{m} d\mu(m) \rangle, \; \; \forall f \in U.
\end{eqnarray*}
Thus, $T^{*}$ is adjointable. So, for any $f\in U$, we obtain
\begin{eqnarray}\label{formula6}
\int_{\mathcal{M}} \langle \Lambda_{m}f,\; \Lambda_{m}f \rangle d\mu(m)=\langle T^{*}f,\; T^{*}f \rangle\leq||T^{*}||^{2} \langle f,\; f \rangle.
\end{eqnarray}
Now, it follows from (\ref{formula5}) and (\ref{formula6})that the family $\{\Lambda_m\}_{m\in \mathcal{M}}$ is a $c$-$K$-$g$-frame for $U$ with respect to $ \{V_{m}\}_{m \in \mathcal{M}}.$
\end{proof}

\begin{thm}\label{thm3}
Suppose that $K\in End^{*}_{\mathcal{A}}(U)$ and
$\{\Lambda_m\}_{m\in \mathcal{M}} \subset End^{*}_{\mathcal{A}}(U, \; V_{m})$ is a $c$-$g$-Bessel
system for $U$ with respect to $\{V_{m}\}_{m \in \mathcal{M}}$  with the synthesis operator $T$ such that $\overline{Ran}(T^{*})$ is orthogonally complemented.
Then, $\{\Lambda_m\}_{m\in \mathcal{M}}$
 is a $c$-$K$-$g$-frame for $U$
 if and only if $R(K)\subseteq R(T)$.
\end{thm}
\begin{proof}
Given a $c$-$g$-orthonormal basis $\Theta=
\{\Theta_{m}\}_{m\in \mathcal{M}}\subset End^{*}_{\mathcal{A}}(U,\;V_{m})$ for $U$ with respect to
$\{V_{m}\}_{m \in \mathcal{M}}$. First, sassume $\{\Lambda_m\}_{m\in \mathcal{M}}$ is a $c$-$K$-$g$-frame for $U$ with respect to
 $\{V_m\}_{m\in \mathcal{M}}$.
 Thus, by definition there exists a positive constant $A$ such that
\begin{eqnarray}\label{formula7}
A \langle K^{*}f, \; K^{*}f \rangle \leq
\int_{\mathcal{M}} \langle \Lambda_{m}f, \;\Lambda_{m}f \rangle d\mu, \; \forall f \in U.
\end{eqnarray}
According to  Lemma \ref{form.22}, we have $\Lambda_{m}=\Theta_{m}T^{*}\; for \; all \; m\in {\mathcal{M}}.$
It follows from (\ref{formula7}) that for any $f \in U$,
\begin{eqnarray*}
A \langle K^{*}f, \; K^{*}f \rangle &\leq&
 \int_{\mathcal{M}} \langle \Theta_{m}T^{*}f, \;\Theta_{m}T^{*}f \rangle d\mu\\
&=&\langle T^{*}f, \; T^{*}f \rangle=\langle TT^{*}f, \; f \rangle,
\end{eqnarray*}
i.e., $AKK ^{*}\leq TT^{*}.$ So , by Lemma \ref{form.2}, we conclude that $R(K)\subseteq R(T)$.\\
Conversely, let $R(K)\subseteq R(T)$. Lemma \ref{form.2}
implies that $KK^{*}\leq \lambda TT^{*}$ for some $\lambda> 0$. Hence, we obtain
\begin{eqnarray}\label{formula8}
\frac{1}{\lambda}\langle K^{*}f, \; K^{*}f \rangle
\leq \langle T^{*}f, \; T^{*}f \rangle, \;\forall f\in U.
\end{eqnarray}
Since $\{\Lambda_m\}_{m\in \mathcal{M}}$ is a $c$-$g$-Bessel system for $U$
and $\{\Theta_m\}_{m\in \mathcal{M}}$ is a $c$-$g$-orthonormal basis for
$U$, using Lemma \ref{form.22}, for every $f \in U$,
 we have
\begin{eqnarray}\label{formula9}
\langle T^{*}f, \; T^{*}f \rangle
= \int_{\mathcal{M}} \langle \Theta_{m}T^{*}f, \;\Theta_{m}T^{*}f \rangle d\mu
=\int_{\mathcal{M}}\langle \Lambda_{m}f, \; \Lambda_{m}f \rangle d\mu.
\end{eqnarray}
Now, from the relations (\ref{formula8}) and (\ref{formula9}), we get
$$\frac{1}{\lambda}\langle K^{*}f, \; K^{*}f \rangle \leq
\int_{\mathcal{M}}\langle \Lambda_{m}f, \; \Lambda_{m}f \rangle d\mu, \; \forall f\in U.$$
 Therefore, $\{\Lambda_{m}\}_{m\in \mathcal{M}}$ is a $c$-$K$-$g$-frame for
 $U$ with respect to $\{V_{m}\}_{m\in \mathcal{M}}.$
\end{proof}

An advantage of studying $c$-$K$-$g$-frames is that a
$c$-$K$-$g$-frame can be constructed using the family of operators which is
not a $c$-$g$-frame:
\begin{exam}
Let $U$ be a Hilbert $C^*$-module, $\mathcal{M} =\mathbb{N}$ and
$\mu$ be the counting measure. Let $\{e_m\}_{m\in \mathcal{M}}$ be
an orthonormal basis for $U$ and $V_m =
\overline{Span}\{e_{3m-2},e_{3m-1},e_{3m}\},\,\ m=1,2,\cdots$. For fixed
$N \in \mathbb{N}$, we define adjointable operator $\Lambda_m
: U \longrightarrow V_m$ as follows:
\begin{eqnarray*}
\Lambda_1 f = \sum_{k=1}^N \langle f, \; e_k \rangle e_k \,\
\mbox{and}\quad\ \Lambda_mf=0 \;\;  for\,\ m\geq 2.
\end{eqnarray*}
Then, it is easy to check  the family $\{\Lambda_m\}_{m \in \mathcal{M}}$ is not a $g$-frame for
$U$ with respect to $\{V_m\}_{m \in \mathcal{M}}$. In fact, if we
take $f=e_{N+1}$, then for every $f\in U$, we get

\begin{eqnarray*}
\sum_{m=1}^{\infty}{\langle\Lambda_mf,\; \Lambda_mf \rangle = \langle
\Lambda_1f,\; \Lambda_1f\rangle = 0}.
\end{eqnarray*}
Define $K:U\longrightarrow U$ by
\begin{eqnarray*}
Ke_m = \left\lbrace
\begin{array}{lr}
me_m, \,\ m\leq N \\
0,\,\ \quad\ m > N.\\
\end{array}
\right.
\end{eqnarray*}
Clearly, $K$ is adjointable and it satisfies:
\begin{center}
$   K^*e_m = \left\lbrace \begin{array}{lr}
 me_m, \,\  m\leq N \\
 0,\,\ \quad\  m > N.\\
 \end{array}
 \right.$

\end{center}
 Now, for any $f \in U$; $ f = \sum_{m=1}^{\infty} {c_m e_m} $, we have
\begin{eqnarray*}
\sum_{m=1}^{\infty}{\langle
\Lambda_{m}f,\; \Lambda_{m}f\rangle}=\langle
\Lambda_{1}f,\; \Lambda_{1}f\rangle =\sum_{m=1}^{N}
c_m c^{*}_m,
\end{eqnarray*}
also
\begin{eqnarray*}
\langle K^* f , \;K^* f \rangle = \sum_{m=1}^{N}{m^2
c_mc^{*}_m}.
\end{eqnarray*}
Therefore, for any $f\in U$, we obtain
\begin{eqnarray*}
\frac{1}{N^2} \langle K^* f, \; K^* f \rangle &=&
\sum_{m=1}^{N}{(\frac{m}{N})^2}  c_mc^{*}_m \leq  \sum_{m=1}^{N}{c_mc^{*}_m}\\
&=&\sum_{m=1}^{\infty}{\langle \Lambda_{m}f,\; \Lambda_{m}f\rangle}
=\langle \Lambda_{1}f,\; \Lambda_{1}f\rangle \\
&\leq& \|\Lambda_1\|^2
\langle f, \;f \rangle.
\end{eqnarray*}
 This shows $\{ \Lambda_{m}\}_{m \in \mathcal{M}}$ is a
$K$-$g$-frame for $U$ with respect to $\{V_m\}_{m \in
\mathcal{M}}$.
\end{exam}
\section{\textbf{Dual $c$-$K$-$g$-Bessel systems}}
In this section, the dual $c$-$K$-$g$-Bessel system of a $c$-$K$-$g$-frame in
Hilbert $C^{*}$-modules is defined. Next, some results on dual $c$-$K$-$g$-Bessel systems are derived. We begin with a result which is the motivation of this definition.
\begin{thm}\label{thm4}
Let $K \in End_{\mathcal{A}}^*(U)$ and $\{\Lambda_m\}_{m\in \mathcal{M}}\subset End_{\mathcal{A}}^*(U, V_m)$
such that for any $f \in U$, the mapping $m \mapsto \langle f, \Lambda_mf \rangle$
is measurable on ${\mathcal{M}}.$  Suppose the operator $T^{*}: U \mapsto \bigoplus_{m\in \mathcal{M}}V_m$  defined by $T^{*}f=\{\Lambda_{m}f\}_{m \in \mathcal{M}}$ and
$\overline{Ran}(T^{*})$ is orthogonally complemented. Then the following statements are equivalent:

\begin{itemize}
\item [(i)] $\{\Lambda_m\}_{m\in \mathcal{M}}$ is a $c$-$K$-$g$-frame for $U$ with respect to $\{V_m\}_{m \in \mathcal{M}}$.
\item [(ii)] $\{\Lambda_m\}_{m\in \mathcal{M}}$ is a $c$-$g$-Bessel system for $U$ with respect to $\{V_m\}_{m \in \mathcal{M}}$ and there exists
a $c$-$g$-Bessel system $\{\Gamma_m\}_{m\in \mathcal{M}}\subset
End_{\mathcal{A}}^*(U, V_m)$  for $U$ with respect to
$\{V_m\}_{m \in \mathcal{M}}$ such that
\begin{eqnarray*}
\langle K f, \;g \rangle=\int_{\mathcal{M}} \langle
\Lambda_m^*\Gamma_m f, \; g\rangle d\mu(m),\,\ \forall f,g \in U.
\end{eqnarray*}
\item [(iii)] The integral $\int_{\mathcal{M}}\Lambda^{*}_{m}g_{m}d\mu(m)$ converges in $U$ for all
$\{g_{m}\}_{m \in \mathcal{M}}\in \bigoplus_{m \in {\mathcal M}}V_m$
and there exists $C> 0$ such that for any $f \in U$, there is $\{g_{m, f}\}_{m \in \mathcal{M}}\in \bigoplus_{m \in {\mathcal M}}V_m$  satisfying:
\begin{eqnarray*}
\langle K f, \; g \rangle=\int_{\mathcal{M}} \langle
\Lambda^{*}_m (g_{m,f}), \; g \rangle d\mu(m), \; \; \forall f,g \in U, \\
and \int_{\mathcal{M}} \langle g_{m,f}, \; g_{m,f} \rangle d\mu(m)\leq C \langle f,\; f \rangle.
\end{eqnarray*}
 \item
\end{itemize}
\end{thm}
\begin{proof}
(i)$\Rightarrow(ii)$. Let $\{\Lambda_m\}_{m\in \mathcal{M}}$ be a
$c$-$K$-$g$-frame for $U$ with respect to $\{V_m\}_{m \in
\mathcal{M}}$ with the bounds $A, B> 0$. Then, for any $ f \in U$,
we get
\begin{eqnarray*}
A\langle K^{*}f, \; K^{*}f \rangle \leq \int_{\mathcal{M}} \langle
\Lambda_m f, \Lambda_m f\; \rangle d\mu(m)\leq B \langle f,\; f \rangle .
\end{eqnarray*}
The left-hand inequality means that $AKK^{*}\leq TT^{*}$. So, by
Lemma \ref{form.2}, there exsist $\Gamma \in
End^{*}_{\mathcal{A}}(U, \bigoplus_{m} V_{m})$ such that
$K=T\Gamma.$ Now, suppose that $P_{n}$ is the projection on
$\bigoplus_{m}V_{m}$ which maps every element to $n$-th
component, i.e., for each $g=\{g_{m}\}_{m \in \mathcal{M}} \in
\bigoplus_{m} V_{m},\; P_{n}g=\{u_{m}\}_{m \in \mathcal{M}}$, where
\begin{eqnarray*}
u_m = \left\lbrace
\begin{array}{lr}
g_n, \,\ m=n \\
0,\,\ \quad\ m \neq n.\\
\end{array}
\right.
\end{eqnarray*}

If we take $\Gamma_{m}=P_{m}\Gamma$, then for any $f \in U$, we
have
\begin{eqnarray*}
\int_{\mathcal M}\langle \Gamma_{m}f, \; \Gamma_{m}f \rangle d\mu(m)&=&
\int_{\mathcal{M}} \langle P_{m} \Gamma f, \; P_{m} \Gamma f \rangle d\mu(m) \\
&=& \int_{\mathcal{M}} \langle (\Gamma f)_{m} , \; (\Gamma f)_{m}  \rangle d\mu(m) \\
&=& \langle \Gamma f, \; \Gamma f \rangle
\leq ||\Gamma||^{2} \langle f,\; f \rangle.
\end{eqnarray*}
Hence, $\{\Gamma_{m}\}_{m \in \mathcal{M}}$ is a $c$-$g$-Bessel system for $U$ with respect to $\{V_{m}\}_{m \in \mathcal{M}}.$
On the other hand, for all $f, g \in U$, we obtain
\begin{eqnarray*}
\langle Kf,\; g \rangle= \langle \Gamma f,\; T^{*}g \rangle &=&
\int_{\mathcal M}\langle P_{m}(\Gamma f), \; \Lambda_{m}g \rangle d\mu(m)\\
&=&\int_{\mathcal{M}} \langle \Gamma_{m} f, \; \Lambda_{m}g \rangle d\mu(m) \\
&=& \int_{\mathcal M}\langle \Lambda_{m}^{*}\Gamma_{m} f, \; g \rangle d\mu(m).
\end{eqnarray*}
$(ii) \Rightarrow (iii)$. Since $\{\Lambda_m\}_{m \in \mathcal{M}}$ is a $c$-$g$-Bessel system for $U$, for every $\{g_m\}_{m \in \mathcal{M}} \in \bigoplus_{m} V_m$, the integral $\int_{\mathcal{M}}\Lambda_{m}^{*}(g_{m}) d\mu(m)$ converges in $U$ with the norm topology. Assume $B$ is the $c$-$g$-Bessel bound of $\{\Gamma_{m}\}_{m \in \mathcal{M}}$ and for any
$f \in U$,  we put $g_{m,f}=\Gamma_{m}f, \; \forall m \in {\mathcal{M}}.$ Then for all $f, g \in U$, we get
\begin{eqnarray*}
\langle Kf,\; g \rangle &=&
\int_{\mathcal{M}} \langle \Lambda_{m}^{*}\Gamma_{m}f, \; g \rangle d\mu(m)\\
&=&\int_{\mathcal{M}} \langle \Lambda_{m}^{*}g_{m, f}, \; g \rangle d\mu(m).
\end{eqnarray*}
 also
\begin{eqnarray*}
\int_{\mathcal M}\langle g_{m,f}, \; g_{m,f} \rangle d\mu(m)
=\int_{\mathcal{M}} \langle \Gamma_{m}f, \; \Gamma_{m}f \rangle d\mu(m)\leq B \langle f,\; f \rangle.
\end{eqnarray*}
$(iii) \Rightarrow (i)$. Since for each $\{g_m\}_{m \in
\mathcal{M}} \in \bigoplus_{m} V_m$, the integral
$\int_{\mathcal{M}}\Lambda_{m}^{*}(g_{m}) d\mu(m)$ converges in
$U$ with the norm topology and $\{\Lambda_{m}f\}_{m \in
\mathcal{M}} \in \bigoplus_{m} V_m, \; \forall f \in U$, we can
conclude as in the proof of Theorem \ref{thm2}, $T^{*}$ is
adjointable. Thus for  every $f \in U$, we obtain
\begin{eqnarray}\label{formula10}
\left|\left| \int_{\mathcal{M}}\langle \Lambda_{m}f,\;\Lambda_{m}f \rangle d\mu(m) \right|\right|=||\langle T^{*}f, \; T^{*}f \rangle||\leq ||T^{*}||^{2} \langle f,\; f \rangle.
\end{eqnarray}
On the other hand, for any $f \in U$,  we get
\begin{eqnarray*}
||K^{*}f||^{2}&=& sup_{g \in U, ||g||=1}||\langle K^{*}f,\; g
\rangle||^{2}\\
&=&sup_{g \in U, ||g||=1}||\langle f,\; Kg \rangle||^{2}\\
&=&sup_{g \in U, ||g||=1} \left|\left| \langle f,\;  \int_{\mathcal{M}}\Lambda_{m}^{*}g_{m,g}d\mu(m) \rangle \right|\right|^{2}\\
&=& sup_{g \in U, ||g||=1} \left|\left| \int_{\mathcal{M}}\langle \Lambda_{m}f,\; g_{m,g} \rangle d\mu(m) \right|\right|^{2}\\
&\leq& sup_{g \in U, ||g||=1} \left|\left| \int_{\mathcal{M}}
\langle \Lambda_{m}f,\; \Lambda_{m}f \rangle d\mu(m) \right|\right| \\
&\times&
\left|\left| \int_{\mathcal{M}} \langle g_{m, g},\;
g_{m, g} \rangle d\mu(m) \right|\right|\\
&\leq& C \left|\left|\int_{\mathcal{M}}\langle \Lambda_{m}f,\; \Lambda_{m}f \rangle
d\mu(m)\right|\right|.
\end{eqnarray*}

Therefore, we have
\begin{eqnarray}\label{formula11}
C^{-1}||K^{*}f||^{2}\leq \left|\left|\int_{\mathcal{M}} \langle \Lambda_{m}f, \; \Lambda_{m}f \rangle d\mu(m) \right|\right|, \; \forall f \in U.
\end{eqnarray}
 Now,  it follows from (\ref{formula10}) and (\ref{formula11}) that $\{\Lambda_{m}\}_{m \in \mathcal{M}}$ satisfies in the condition of Theorem \ref{thm2}.
 Hence, $\{\Lambda_m\}_{m \in \mathcal{M}}$ is a $c$-$K$-$g$-frame for $U$ with respect to $\{V_{m}\}_{m \in \mathcal{M}}.$
\end{proof}

\begin{defn}
Let $K \in End_{\mathcal{A}}^*(U)$ and $\{\Lambda_m\}_{m \in
\mathcal{M}} \subset End^{*}_{\mathcal{A}}(U,\; V_{m})$ be a $c$-$K$-$g$-frame for $U$ with respect to $\{V_m\}_{m \in \mathcal{M}}$.
 A $c$-$g$-Bessel system
$\{\Gamma_m\}_{m \in \mathcal{M}}\subset End_{\mathcal{A}}^*(U,
V_m)$ is called a dual $c$-$K$-$g$-Bessel system of $\{\Lambda_m\}_{m \in
\mathcal{M}}$ if it satisfies:
\begin{eqnarray*}
\langle K f, \;g \rangle=\int_{\mathcal{M}} \langle
\Lambda_m^*\Gamma_m f, \;g\rangle d\mu(m),\; \; \forall f,g \in U.
\end{eqnarray*}
In this case,  we say that $(\Lambda_m, \; \Gamma_m)$ forms a
$c$-$K$-$g$-dual pair for $U$ with respect to $\{V_m\}_{m \in \mathcal{M}}$.
\end{defn}
\begin{prop}
Let $K \in End_{\mathcal{A}}^*(U)$ and $\{\Lambda_m\}_{m \in
\mathcal{M}} \subset End^{*}_{\mathcal{A}}(U, V_{m})$ be a $c$-$K$-$g$-frame for $U$ with respect to
$\{V_m\}_{m \in \mathcal{M}}$  with the $c$-$g$-frame operator $S$.
 Then, for each
$\alpha\in \mathbb{R}$, $(\Lambda_m S^{\alpha}, \;\Lambda_m
S^{-1-\alpha}K)$ forms a $c$-$K$-$g$-dual pair for $U$ with respect to
$\{V_m\}_{m \in \mathcal{M}}$.
\end{prop}
\begin{proof}

Firstly, assume  $\{\Lambda_m\}_{m \in \mathcal{M}}$ is a $c$-$K$-$g$-frame for $U$ with the upper bound $B$.
Thus, for any $\alpha \in {\mathbb{R}}$ and any $f \in U$, we have

$$\int_{\mathcal{M}} \langle \Lambda_mS^{\alpha}f, \; \Lambda_m S^{\alpha}f\rangle
d\mu(m)\leq B \langle S^{\alpha}f, \;S^{\alpha}f \rangle
\leq B^{2\alpha+1}\langle f, \;f \rangle.$$
Hence, the family $\{\Lambda_mS^{\alpha}\}_{m \in \mathcal{M}}$
is a c-$g$-Bessel system for $U$ with respect to
$\{V_m\}_{m \in \mathcal{M}}$.
On the other hand, for all $f, g\in U$, we get
\begin{eqnarray*}
\langle K f, \; g \rangle &=& \langle S^{\alpha}S S^{-1-\alpha}K f, \;g
\rangle\\ &=&\langle SS^{-1-\alpha}K f, \; S^{\alpha}g
\rangle\\
&=&\int_{\mathcal{M}}\langle \Lambda^*_m\Lambda_m S^{-1-\alpha}K
f, \; S^{\alpha}g \rangle d\mu(m)\\&=&\int_{\mathcal{M}}\langle
\Lambda_m S^{-1-\alpha}K f, \; \Lambda_m S^{\alpha}g \rangle
d\mu(m).
\end{eqnarray*}
Therefore, for all $\alpha\in \mathbb{R}$, $(\Lambda_mS^{\alpha}, \; \Lambda_mS^{-1-\alpha}K)$ forms a
$c$-$K$-$g$-dual pair for $U$ with respect to
$\{V_m\}_{m \in \mathcal{M}}$. In particular, the family
$(\Lambda_m, \; \Lambda_mS^{-1})$ forms a $c$-$K$-$g$-dual pair for $U$ with respect to
$\{V_m\}_{m \in \mathcal{M}}$.
\end{proof}
The family $\{\Lambda_mS^{-1}\}_{ m \in \mathcal{M}}$ is called
$c$-$K$-$g$-canonical dual of $\{\Lambda_m\}_{ m \in \mathcal{M}}$.

The next result shows that the property of dual $c$-$K$-$g$-Bessel system preserves under a co-isometry.
\begin{thm}
Let $\{\Lambda_m\}_{m \in \mathcal{M}}$ and $\{\Gamma_m\}_{m \in \mathcal{M}}
\subset End^{*}_{\mathcal {A}}(U, V_{m})$  such that $(\Lambda_{m},\; \Gamma_{m})$
forms a $c$-$K$-$g$-dual pair for $U$ with respect to $\{V_{m}\}_{m \in \mathcal{M}}$.
If  $\Theta \in End^{*}_{\mathcal{A}}(U)$ is a co-isometry with $K\Theta
=\Theta K$, then $(\Lambda_{m}\Theta^{*},\; \Gamma_{m}\Theta^{*})$ forms a
$c$-$K$-$g$-dual pair for $U$ with respect to $\{V_{m}\}_{m \in \mathcal{M}}$.
\end{thm}
\begin{proof}
Since $\Theta$ is a co-isometry,  we have
 $ \langle \Theta^{*}f,\; \Theta^{*}g \rangle=
 \langle f, \; g \rangle,\; \forall f, g \in U, \; or \; \Theta\Theta^{*}=I_{U}$.
 By  definition there exist positive numbers $A, B$ such that
$$A \langle K^{*}f,\; K^{*}f \rangle
\leq \int_{\mathcal{M}} \langle \Lambda_{m}f,\; \Lambda_{m}f \rangle d\mu(m)
\leq B \langle f,\; f \rangle, \; \; \forall f \in U.$$
Thus, for any $f \in U$, we get
$$\int_{\mathcal{M}} \langle \Lambda_{m}\Theta^{*}f,\; \Lambda_{m}\Theta^{*}f \rangle d\mu(m)
\leq B \langle \Theta^{*}f,\; \Theta^{*}f \rangle=B\langle f, \; f \rangle.$$
On the other hand, for any $f \in U$, we have
\begin{eqnarray*}
\int_{\mathcal{M}} \langle \Lambda_{m}\Theta^{*}f,\; \Lambda_{m}\Theta^{*}f \rangle d\mu(m)
&\geq& A \langle K^{*}\Theta^{*}f,\; K^{*}\Theta^{*}f \rangle \\
&=& A \langle \Theta^{*}K^{*}f,\; \Theta^{*}K^{*}f \rangle \\
&=&A \langle K^{*}f,\; K^{*}f \rangle.
\end{eqnarray*}

Hence, $\{\Lambda_{m}\Theta^{*}\}_{m \in \mathcal{M}}$ is a $c$-$K$-$g$-frame for $U$
with respect to $\{V_{m}\}_{m \in \mathcal{M}}$.
Also, since $\{\Gamma_m\}_{m \in \mathcal{M}}$ is a $c$-$K$-$g$-Bessel system for $U$,
for a positive constant $C$, we get
$$ \int_{\mathcal{M}} \langle \Gamma_{m}f,\; \Gamma_{m}f \rangle d\mu(m)
\leq C \langle f,\;f \rangle, \; \forall f \in U.$$
Thus, for any $f \in U$,  we obtain
$$ \int_{\mathcal{M}} \langle \Gamma_{m}\Theta^{*}f,\; \Gamma_{m}\Theta^{*}f \rangle d\mu(m)
\leq C \langle \Theta^{*}f,\; \Theta^{*}f \rangle=C \langle f,\; f \rangle.$$

Hence, $\{\Gamma_{m}\Theta^{*}\}_{m \in \mathcal{M}}$ is a $c$-$g$-Bessel
system for $U$ with respect to $\{V_{m}\}_{m \in \mathcal{M}}$. Because
$(\Lambda_{m},\; \Gamma_{m})$ forms a  $c$-$K$-$g$-dual pair for $U$,  we have
$$ \langle Kf,\; g \rangle
=\int_{\mathcal{M}} \langle \Lambda_{m}^{*}\Gamma_{m}f,\; g \rangle d\mu(m), \; \; \forall f, g \in U.$$

Therefore, for all $f, g \in U$, we get

\begin{eqnarray*}
\int_{\mathcal{M}} \langle (\Lambda_{m}\Theta^{*})^{*}(\Gamma_{m}\Theta^{*})f,\; g\rangle d\mu(m)
&=&\int_{\mathcal{M}} \langle \Lambda_{m}^{*}\Gamma_{m}\Theta^{*}f,\; \Theta^{*}g \rangle d\mu(m) \\
&=& \langle K\Theta^{*}f,\; \Theta^{*}g \rangle
=\langle \Theta K\Theta^{*}f,\; g\rangle \\
&=&\langle K\Theta\Theta^{*}f,\; g \rangle=\langle Kf, \; g \rangle.
\end{eqnarray*}
So, $(\Lambda_{m}\Theta^{*}, \Gamma_{m}\Theta^{*})$ forms a $c$-$K$-$g$-dual
pair for $U$ with respect to $\{V_{m}\}_{m \in \mathcal{M}}$.
\end{proof}
\begin{thm}\label{thm5}
Let $K\in End^{*}_{\mathcal{A}}(U)$ and
 $\{\Theta_m\}_{m\in \mathcal{M}} \subset End^{*}_{\mathcal{A}}(U,\; V_{m})$ be a $c$-$g$-orthonormal basis
 for $U$ with respect to $\{V_{m}\}_{m\in \mathcal{M}}$. Suppose that
$\{\Lambda_{m}\}_{m\in \mathcal{M}} \subset End^{*}_{\mathcal{A}}(U,\; V_{m})$ is a $c$-$K$-$g$-frame for $U$ and $\{\Gamma_{m}\}_{m\in \mathcal{M}}\subset End^{*}_{\mathcal{A}}(U, \; V_{m}) $ is a $c$-$g$-Bessel system for $U$ with respect to $\{V_{m}\}_{m\in \mathcal{M}}$ with the synthesis operators $T$ and $R$, respectively.
Then, $\{\Gamma_{m}\}_{m\in \mathcal{M}}$
is a  dual $c$-$K$-$g$-Bessel system of $\{\Lambda_{m}\}_{m\in \mathcal{M}}$ if and only if $K=TR^{*}.$
\end{thm}
\begin{proof}
Assume that $\{\Gamma_{m}\}_{m\in \mathcal{M}}$ is a dual $c$-$K$-$g$-Bessel system
of $\{\Lambda_{m}\}_{m\in \mathcal{M}}$.
Using Lemma \ref{form.22}, we have $\Lambda_{m}=\Theta_{m}T^{*}$ and
$\Gamma_{m}=\Theta_{m}R^{*}, \; \forall m\in \mathcal{M}$.
Now, by the definition, for all $f, g \in U$, we get
\begin{eqnarray*}
\langle Kf,\; g \rangle=\int_{\mathcal{M}}\langle \Lambda^{*}_{m}\Gamma_{m}f,\; g \rangle d\mu
&=&\int_{\mathcal{M}} \langle (\Theta_{m}T^{*})^{*}(\Theta_{m}R^{*})f,\; g \rangle d\mu \\
&=&\int_{\mathcal{M}}\langle \Theta_{m}R^{*}f,\; \Theta_{m}T^{*}g \rangle d\mu \\
&=&\langle R^{*}f,\; T^{*}g \rangle
=\langle TR^{*}f, \; g \rangle.
\end{eqnarray*}
It follows from  that $K=TR^{*}.$\\
Conversely, suppose  $K=TR^{*}$. Because $\{\Gamma_{m}\}_{m\in \mathcal{M}}$ and
 $\{\Lambda_{m}\}_{m\in \mathcal{M}}$ are $c$-$g$-Bessel systems for $U$. Again,
 by Lemma \ref{form.22}, we have
 $\Lambda_{m}=\Theta_{m}T^{*}$ and $\Gamma_{m}=\Theta_{m}R^{*} \; for\; all\; m\in \mathcal{M}.$

 Thus, similar to the above argument, we obtain
 $$ \int_{\mathcal{M}} \langle \Lambda^{*}_{m}\Gamma_{m}f, \; g \rangle= \langle TR^{*}f,\; g \rangle
=\langle Kf, \; g \rangle, \; \forall f, g\in U.$$  Therefore, $\{\Gamma_{m}\}_{m\in \mathcal{M}}$
 is a  dual $c$-$K$-$g$-Bessel system of $\{\Lambda_{m}\}_{m\in \mathcal{M}}.$
\end{proof}
{\bf Declaration}: All authors declare that they have no conflicts of interest.

\end{document}